%% file: No-Speeduparxiv.tex
\documentclass[a4paper,12pt]{article}
\usepackage{amssymb}
\usepackage{latexsym}
\usepackage{prooft}
\usepackage{bussproofs}
\usepackage{color}
\input{Barr_Remarks_macros}

\newcommand{\citep}{\cite}

\newcommand{\Em}{\mathsf{E}}
\newcommand{\nege}{\neg_{\Em}}
\newcommand{\fa}{\varphi}
\newcommand{\fb}{\psi}
\newcommand{\fc}{\theta}
\newcommand{\IQ}{{\mathcal Q}}
\newcommand{\IR}{{\mathcal R}}
\newcommand{\IJJ}{{\mathcal J}}
\newcommand{\uml}{\"}

\begin{document}

\title{No speedup for geometric theories\footnote{In honor of Grisha Mints.}}

\author{Michael Rathjen\\
{\small Department of Pure Mathematics, University of Leeds}\\
{\small Leeds LS2 9JT, United Kingdom}\\
 {\tt \scriptsize
M.Rathjen@leeds.ac.uk}}

\maketitle
\begin{abstract} Geometric theories based on classical logic are conservative over their intuitionistic counterparts for geometric implications. The latter result (sometimes referred to as Barr's theorem) is squarely a consequence of Gentzen's Hauptsatz. Prima facie though, cut elimination can result in superexponentially longer proofs. In this paper it is shown that the transformation 
of a classical proof of a geometric implication in a geometric theory into an intuitionistic proof can be achieved  in feasibly many steps.

\end{abstract}

\section{Introduction} For  geometric theories it is known that the existence of a classical proof  of a geometric implication yields the existence of an intuitionistic proof. Existing effective proofs of this fact  use cut elimination (see \cite{orevkov,negplat2011,rathjen2016}) and thus are liable to  effect a superexponential blow-up. When I visited Stanford in November 2013, Grisha Mints told me 
 that he was aiming to show that 
the transformation can be achieved via a polynomial time algorithm. Sadly, this was the last time I saw him. Some of what he told me is reflected in his  article \cite{mints}, which presents a partial result in that a merely polynomially longer proof is achievable in an augmentation of the intutitionistic theory via relational Skolem axioms.\footnote{See also the comments on \cite{mints} in \cite{dyckneg}.} The article, though, does not achieve what he had originally intended. At the time, I did not start to work on the problem  but the conversation with
Grisha lingered in my mind over the years  and for some reason I became convinced that there should be a simple argument based on an interpretation akin to the Friedman-Dragalin $A$-translation.
The short paper before the reader is an elaboration of that idea. 

As I was working on such a translation it dawned on me that similar ideas must have been considered before.
Indeed, Leivant in \cite{leivant85} used interpretations of classical logic into intuitionistic logic to obtain (partial) conservativity results for classical theories over their intuitionistic 
versions. He explained that the Friedman-Dragalin interpretation ``is derived naturally from a trivial translation of $I$ into $M$'' (\cite[p. 683]{leivant85}), where $I$ and $M$ signify intuitionistic and minimal logic, respectively,
and ``trivial translation'' refers to Kolmogorov's 1925 translation \cite{kolmogorov25}. 
The observation that the latter is actually an interpretation into Johanssons' 1937 minimal logic \cite{johansson37} emerges as the fulcrum for obtaining conservativity results in that in minimal
logic falsum $\perp$ just acts as placeholder for an arbitrary formula. Crucial in the machinery of \cite{leivant85} are the definitions of three syntactical classes, namely  the {\em spreading, wiping and isolating} schemata and formulas with regard to the Kolmogorov interpretation (also see \cite{TvD88}, Ch.2, Sect. 3 for an exposition).  
The ideas underlying these classes have informed Definition \ref{I6} and Proposition \ref{I7}.  
But alas I couldn't see how to directly infer the conservativity of classical over intuitionistic geometric theories by reassembling results from \cite{leivant85}. 
However, Ishihara's article \cite{ishihara2000} actually furnishes what is needed. A crucial move in \cite{ishihara2000} is to introduce a new propositional constant into  the language which can act as a placeholder for arbitrary formulas not only in minimal logic but also in intuitionistic logic.

\section{Geometric theories}
\begin{deff}\label{geometricp}{\em The {\em positive formulas} are 
constructed from atomic formulas and falsum $\perp$ by $\wedge$, $\vee$, and $\exists$.

{\em Geometric implications} are made up of the positive formulas, implications of positive formulas and the result of prefixing universal quantifiers to positive formulas and implications of positive formulas. 

$\neg \varphi$ is defined as $\varphi \to \perp$. Thus if $\varphi$ is a positive formula then $\neg \varphi$ is a geometric  implication.

 A theory is {\em geometric} if all its axioms are geometric implications.
 }\end{deff}

Below we shall give several examples of geometric theories. 
 \begin{examples}\label{3.11p}{\em
 \begin{itemize}
 \item[(i)] 1. Robinson arithmetic formulated in the language with a constant $0$, a unary successor
function symbol $\suc$, binary function symbols $+$ and $\cdot$, and a binary predicate symbol $<$. 
\comment{Axioms are the equality axioms and the universal closures of the following.
\begin{enumerate}
\item $\neg \suc(a) = 0$.
\item $\suc(a) = \suc(b)\to a = b$.
\item $a = 0\,\vee\,\exists y \,a = \suc(y)$.
\item $ a + 0 = a$.
\item $ a + \suc(b) = \suc(a + b)$.
\item $a \cdot 0 = 0$.
\item $ a \cdot \suc(b) = a \cdot b + a$
\end{enumerate}
A classically equivalent axiomatization is obtained if (3) is replaced by
$$ \neg a = 0 \to \exists y\, a = \suc(y)$$ but this is not a geometric implication.}

\item[(ii)] The theories of  groups, rings, local rings and division rings
have geometric axiomatizations. Local rings are commutative rings with $0\ne 1$ having just one maximal ideal. On the face of it,
the latter property appears to be second order but it can be rendered geometrically as follows:
$$\forall x\,(\exists y\,x\cdot y=1\;\vee\;\exists y\,(1-x)\cdot y=1).$$

\item[(iii)] The theories of fields, ordered fields, algebraically closed fields and
real closed fields have geometric axiomatizations.
To express invertibility of non-zero elements one uses $\forall x\,(x=0\,\vee\,\exists y\,x\cdot y=1)$ rather
than the non-geometric axiom $\forall x\,(x\ne 0\to \exists y\,x\cdot y=1)$.

To express algebraic closure replace axioms

$$s\ne 0\,\to \,\exists x\,sx^n+t_1x^{n-1}+\ldots +t_{n-1}x+t_n=0$$
by
$$s= 0\,\vee \,\exists x\,sx^n+t_1x^{n-1}+\ldots +t_{n-1}x+t_n=0$$
where $sx^k$ is short for $s\cdot x\cdot\ldots \cdot x$ with $k$ many $x$.
\\[1ex]
 Also the theory of {\em differential fields} has a geometric axiomatization.
This theory is written in the language of rings with an additional unary function symbol $\delta$. The axioms
are the field axioms plus $\forall x\forall y\,\delta (x+y)=\delta (x)+\delta (y)$ and
$\forall x\forall y\,\delta (x\cdot y)=x\cdot\delta (y)+y\cdot\delta (x)$.

\item[(iv)] The theory of projective geometry has a geometric axiomatization.

\item[(v)] The theories of equivalence relations, dense linear orders, infinite sets
and graphs also have geometric axiomatizations.

\end{itemize}

 }\end{examples}
 
 It is also interesting that Kant's logic in  his {\em Critique of Pure Reason} \cite{kant1} and the {\em J\"asche Logik} \cite{kant2} can be identified with geometric logic 
 as shown by T.~Achourioti and M. van Lambalgen in \cite{lambalgen,lambalgen2}.

 \comment{Kant's logic is to  to be founded on the act of
judging and the different forms of judgement, hence, take pride of place in his
argumentation. The consensus view is that this aspect of the Critique of Pure
Reason is a failure because Kant’s logic is far too weak to bear such a weight.
Here we show that the consensus view is mistaken and that Kant’s logic should
be identified with geometric logic, a fragment of intuitionistic logic of great
foundational significance.}

\section{Conservativity}

The best and most elegant proof system for proof-theoretic investigations  is Gentzen's sequent calculus. With minor notational variations, this article will follow the presentation in Takeuti's  book {\em Proof Theory}
\cite{takeuti87}.  We will deviate, though,  a bit from the setup in chapter 1 of \cite{takeuti87} in that we 
\begin{itemize} 
\item use $\perp$ as a propositional constant (or 0-ary predicate symbol) and define $\neg \varphi$ to be $\varphi \to \perp$;
\item use the symbol $\to$ rather than $\supset$ for the implication symbol;
\item use $\Rightarrow$ to separate the left and right part of a sequent, i.e., $\Gamma\Rightarrow  \Delta$ rather than $\Gamma\to \Delta$;
\item add  sequents $\Gamma,\perp \Rightarrow \Delta$ as  axioms\footnote{In \cite{takeuti87}, axioms are called {\em initial sequents}.} and omit the rules for $\neg$ (this axiom scheme for $\perp$ will be refereed to as  $\mathsf{Ax}_{\perp}$).
\end{itemize}
\comment{The formula formation rules in \cite{takeuti87} forbid overlapping scopes of quantifiers with the same variable. This happens for instance in $\forall x\,(P(x)\,\wedge\,\exists x\,Q(x))$. While there is no real loss involved in outlawing such formulas ($\forall y\,(P(y)\,\wedge\,\exists x\,Q(x))$ semantically fulfills the same job), this requirement is a bit awkward 
when one wants to replace a placeholder propositional constant by a formula. As a result, we allow overlapping scopes of quantifiers with the same quantified variable. }

\begin{deff}{\em 
Intuitionistic sequents $\Gamma\Rightarrow \Delta$ satisfy the extra requirement that the succedent $\Delta$ contains at most one formula. In the {\bf intuitionistic} version of this sequent calculus 
only intuitionistic sequents are allowed. In the {\bf minimal logic} version only intuitionistic sequents are permitted  and the scheme $\mathsf{Ax}_{\perp}$ is omitted.

We convey derivability of a sequent $\Gamma\Rightarrow\Delta$ in classical, intuitionistic, and minimal logic by writing $\vdash_c \Gamma\Rightarrow\Delta$,
$\vdash_i \Gamma\Rightarrow\Delta$, and $\vdash_m \Gamma\Rightarrow\Delta$, respectively.

A theory $T$ is a set of sentences. In derivability in $T$ one can use any sequent $\Gamma\Rightarrow \varphi$ with $\varphi\in T$ as  an axiom  (initial sequent). 
 $T\vdash_c \Gamma\Rightarrow\Delta$,
$T\vdash_i \Gamma\Rightarrow\Delta$, and $T\vdash_m \Gamma\Rightarrow\Delta$ are defined accordingly.

For a formula $\varphi$, we shall write $\vdash_c \varphi$,
$\vdash_i \varphi$, and $\vdash_m \varphi$ to convey that  $\vdash_c \emptyset\Rightarrow\varphi$, $\vdash_c \emptyset\Rightarrow\varphi$, and $\vdash_c \emptyset\Rightarrow\varphi$, respectively, where $\emptyset$ stands for the empty sequence of formulas. 
}
\end{deff}

 \begin{deff}\label{I3}{\em 
      We shall use $\Em$ as a symbol for a new propositional constant (or predicate symbol of arity 0). Its purpose will be to serve as a placeholder for an arbitrary formula. Let $\nege \varphi$ be an abbreviation for $\varphi \to \Em$. 
      The $\Em$-negative translation $^\Em$ is defined as follows:
      \begin{eqnarray*} P^\Em:= \nege\nege P\mbox{ for $P$ prime and $P\not\equiv \perp$}; && \perp^\Em := \Em \\
      (\varphi \circ \psi)^\Em := \varphi^\Em\circ \psi^\Em\mbox{ for }\circ \in \{\wedge,\to\}; &&   (\varphi \vee\psi)^\Em := \nege\nege (\varphi^\Em\vee \psi^\Em) \\
      (\forall x \,\varphi)^\Em := \forall x\,\varphi^\Em; &&  (\exists x\,\varphi)^\Em :=\nege\nege \exists x\,\varphi^\Em.
      \end{eqnarray*}
  }\end{deff}
  The foregoing translation is  basically the Gentzen-G\uml odel negative translation (see \cite[3.4, 3.5]{TvD88}), which engineers an interpretation of classical logic into minimal logic.
 
  \comment{
  \begin{prop}\label{I4} \begin{itemize}
  \item[(i)] $\vdash_m \neg\neg A^{\Em}\leftrightarrow A^g$.
  \item[(ii)] $\Gamma\vdash_c\fa \;\Rightarrow\;\Gamma^g\vdash_m\fa^g$.
  \end{itemize}
  \end{prop}
  \prf Since $vdash_m \nege\nege(\theta_0\vee\theta_1)\leftrightarrow \nege (\nege \theta_0\,\wedge\,\nege \theta_1)$ and $\vdash_m \nege\nege\exists x\theta(x)\leftrightarrow \nege \forall x\nege\theta(x)$  it follows that the $^\Em$ translation amounts to the same as Gentzen-G\uml odel $^g$ translation (see \cite[3.4, 3.5]{TvD88}), taking into account 
  
  \cite[Theorem 2.3.5]{TvD88}. \qed
}

  \begin{cor}\label{I5} Given a theory $T$, the theory $T^{\Em}$ has as axioms all formulas $\psi^{\Em}$ with $\psi$ an axiom of $T$.
  \begin{itemize}
  \item[(i)] $\vdash_m\nege\nege \fa^\Em\leftrightarrow \fa^\Em$.
  \item[(i i)] $T\vdash_c\fa \;\Rightarrow\;T^\Em\vdash_m\fa^\Em$.
  \end{itemize}
  \end{cor}
  \prf   Since $\vdash_m \nege\nege(\theta_0\vee\theta_1)\leftrightarrow \nege (\nege \theta_0\,\wedge\,\nege \theta_1)$ and $\vdash_m \nege\nege\exists x\theta(x)\leftrightarrow \nege \forall x\nege\theta(x)$ hold it follows that  the $^\Em$-translation amounts to the same as Gentzen-G\uml odel-$^g$ translation (see \cite[3.4, 3.5]{TvD88}), taking into account 
   that in minimal logic $\perp$ is an arbitrary propositional constant for which we can substitute $\Em$.
  Hence (i) follows from \cite[2.3.3]{TvD88} and (ii) from \cite[2.3.5]{TvD88}. \qed \comment{and (Proposition \ref{I4}(i) and (ii) from \ref{I4}(ii). }

  Below we shall frequently adopt the convention that a string of implications $\to$  is considered to be bracketed to the right, i.e., 
  $\varphi_1\to \varphi_2\to \ldots \to\varphi_{n-1} \to \varphi_n$ is an abbreviation for $\varphi_1\to(\varphi_2\to(\ldots(\varphi_{n-1}\to\varphi_n)\ldots))$
  
  \begin{lem}\label{I2}
  \begin{enumerate}
  \item $\vdash_m \fa\to \nege \nege \fa$;
  \item $\vdash_m(\fa\to\fb)\to(\nege\nege \fa\to \nege\nege \fb)$;
  \item $\vdash_m(\nege\neg (\fa\wedge\fb) \to \nege\nege (\fa\wedge\fb)$;
  \item $\vdash_m\nege\nege \fa\,\wedge\,\nege\nege \fb \to \nege\nege (\fa\wedge\fb)$;
  \item $\vdash_m\nege\neg(\fa\vee\fb)\to \nege\nege(\nege \neg\fa\,\vee\,\nege\neg\fb)$;
  \item $\vdash_m \nege\nege(\nege\nege\fa\,\vee\,\nege\nege\fb)\to \nege\nege(\fa\vee\fb)$;
  \item $\vdash_m\nege\neg(\fa\to\fb)\to (\nege\nege\fa\to\nege\neg\fb)$;
  \item $\vdash_i(\nege \neg \fa\to\nege\nege\fb)\to \nege\nege (\fa\to\fb)$;
  \item $\vdash_m\nege \neg\forall x\fa(x)\to \forall x\nege\neg\fa(x)$;
  \item $\vdash_m\nege\nege \exists x\nege \nege\fa(x)\to \nege\nege\exists x\fa(x)$.
  \end{enumerate}
  \end{lem}
  \prf These claims are stated in \cite[Lemma 2]{ishihara2000} without  proofs.
   (1), (2), (4), (6), and (10) are wellknown with $\Em$ replaced by $\perp$ (see e.g. \cite[1.2]{leivant85}), so it's clear that they hold in minimal logic.
   We now turn to the interesting cases that mix $\neg$ and $\nege$.
  \\[1ex]
  For (3), notice that $\vdash_m \neg \varphi \to \neg(\varphi\wedge \psi)$ and $\vdash_m \neg \psi \to \neg(\varphi\wedge \psi)$, and therefore
  $$\vdash_m(\neg(\varphi\wedge\psi)\to \Em)\to(\neg \varphi\to\Em)\,\wedge\,(\neg \psi\to \Em).$$
  
  (5): We have $\vdash_m[(\varphi\to\perp)\,\vee\,(\psi\to\perp)]\to (\varphi\vee\psi)\to\perp$, yielding
  $$\vdash_m [((\varphi\vee\psi)\to\perp)\to\Em] \to ([(\varphi\to\perp)\,\vee\,(\psi\to \perp)]\to\Em),\mbox{ hence}$$
  \begin{eqnarray*} && \vdash_m[((\varphi\vee\psi)\to\perp)\to\Em]\to((\varphi\to\perp)\to\Em) \\
   && \vdash_m[((\varphi\vee\psi)\to\perp)\to\Em]\to((\psi\to\perp)\to\Em)\end{eqnarray*}
   and thus
   $$\vdash_m [((\varphi\vee\psi)\to\perp)\to\Em]\to [([((\varphi\to\perp)\to\Em)\,\vee\,((\psi \to\perp)\to\Em)]\to\Em)\to\Em].$$
   (7): Successively we see that: 
   \begin{eqnarray*} && \vdash_m\neg\psi\to (\varphi \to \neg(\varphi\to\psi)) \\
    && \vdash_m\neg\psi\to \varphi \to (\neg(\varphi\to\psi)\to\Em) \to\Em\\
    && \vdash_m\neg\psi \to (\neg(\varphi\to\psi)\to\Em) \to\varphi\to\Em\\
     && \vdash_m\neg\psi \to (\neg(\varphi\to\psi)\to\Em) \to((\varphi\to\Em)\to\Em)\to\Em\\
     && \vdash_m (\neg(\varphi\to\psi)\to\Em) \to((\varphi\to\Em)\to\Em)\to\neg\psi\to \Em
     \end{eqnarray*}
     
     (8): $\vdash_i \neg \varphi\to \varphi\to\psi$ and $\vdash_i((\varphi\to\psi)\to\Em)\to\neg\varphi\to\Em$, so
     \begin{eqnarray*} (a) && \vdash_i((\varphi \to\psi)\to\Em)\to[(\neg\varphi\to\Em)\to (\psi\to\Em)\to\psi]\to (\psi\to\Em)\to\Em \\
     (b) && \vdash_i ((\varphi\to \psi)\to\Em) \to \psi\to\Em 
     \end{eqnarray*}
     From (a) and (b) we obtain the desired
     $$\vdash_i((\varphi \to\psi)\to\Em)\to[(\neg\varphi\to\Em)\to (\psi\to\Em)\to\psi]\to\Em.$$
     (9): We have $\vdash_m(\varphi(a)\to\perp)\to \forall x\varphi(x)\to \perp$, and hence 
     \begin{eqnarray*} &&\vdash_m[(\forall x \varphi(x) \to\perp) \to\Em]\to (\varphi(a)\to\perp) \to\Em,\;\;\;\mbox{ whence}\\
     && \vdash_m[(\forall x \varphi(x) \to\perp) \to\Em]\to \forall x [(\varphi(x)\to\perp) \to\Em].
     \end{eqnarray*}
    \qed 

  The following syntactic classes bear some resemblance to   the {\em spreading, wiping and isolating} schemata and formulas in \cite{leivant85} but are actually singled out in 
  \cite{ishihara2000}.

  \begin{deff}\label{I6}{\em We define syntactic classes of formulas $\IQ$, $\IR$ and $\IJJ$ simultaneously by the following clauses:
  \begin{enumerate} \item $\perp$ and every atomic formula $Pt_1\ldots t_n$ belong to $\IQ$. If $Q,Q',\tilde{Q}(a)\in \IQ$ then so are $Q\wedge Q'$, $Q\vee Q'$, $\exists x \tilde{Q}(x)$ and $\forall x \tilde{Q}(x)$. 
  If $Q\in \IQ$ and $J\in \IJJ$ then $J\to Q\in \IQ$.
   \item $\perp\in \IR$.  If $R,R',\tilde{R}(a)\in \IR$ then so are $R\wedge R'$, $R\vee R'$ and $\forall x \tilde{R}(x)$. 
  If $R\in \IR$ and $J\in \IJJ$ then $J\to R\in \IR$.
  
   \item $\perp$ and every atomic formula $Pt_1\ldots t_n$ belong to $\IJJ$. If $J,J',\tilde{J}\in \IJJ$ then so are $J\wedge J'$, $J\vee J'$  and $\exists x \tilde{J}(x)$. 
  If $J\in \IJJ$ and $R\in \IR$ then $R\to J\in \IJJ$.
  \end{enumerate} 
    }\end{deff}
    
    \begin{cor}\label{3.7} \begin{itemize}
    \item[(i)]  All positive formulas are in both, $\IQ$ and $\IJJ$.
    \item[(ii)] All geometric implications are in  $\IQ$.
    \end{itemize}
    \end{cor}
    \prf Obvious. \qed
    
    The following proposition is due to Ishihara \cite{ishihara2000}. 
    
    \begin{prop}\label{I7} \begin{itemize}
    \item[(i)] For $\fa\in \IQ$, $\vdash_i\fa\to \fa^\Em$.
    
     \item[(ii)]  For $\fb\in \IR$, $\vdash_i\nege \neg\fb\to \fb^\Em$.
      \item[(iii)]  For $\fc\in \IJJ$, $\vdash_i \fc^\Em\to \nege \nege\fc$.
  \end{itemize}
  \end{prop} 
  \prf We prove these derivabilities simultaneously by induction on the generation of the classes $\IQ,\IR,\IJJ$. The proof given here is more detailed than the one for
  \cite[Proposition 7]{ishihara2000}. 
  \\[1ex]
  (i): Obviously we have $\vdash_i\perp\to \Em$ and $\vdash_i \psi\to(\psi\to\Em)\to\Em$, which yields $\vdash_i A\to A^{\Em}$ for atomic formulas $A$.
  
  Now suppose $\vdash_i Q_i\to Q_i^{\Em}$ for $i\in\{0,1\}$. Then $\vdash_iQ_0\wedge Q_1\to Q_0^{\Em}\wedge Q_1^{\Em}$, thus $\vdash_iQ_0\wedge Q_1\ to (Q_0\wedge Q_1)^{\Em}$. Likewise one has $\vdash_iQ_0\vee Q_1\to Q_0^{\Em}\vee  Q_1^{\Em}$ and hence $\vdash_iQ_0\vee Q_1\to \nege\nege (Q_0^{\Em}\vee  Q_1^{\Em})$, i.e.,
  $\vdash_iQ_0\vee Q_1\to (Q_0\vee  Q_1)^{\Em}$.
  
  Next assume $\vdash_i Q(a)\to Q(a)^{\Em}$. Then $\vdash_i\forall x Q(x) \to\forall x Q(x)^{\Em}$, so \\ $\vdash_i\forall x Q(x) \to(\forall x Q(x))^{\Em}$. Likewise we have
  $\vdash_i\exists x Q(x) \to\exists x Q(x)^{\Em}$, and so $\vdash_i\exists x Q(x) \to \nege\nege \exists x Q(x)^{\Em}$, which is
  $\vdash_i\exists x Q(x) \to (\exists x Q(x))^{\Em}$.

  Finally assume $\vdash_i J^{\Em}\to \nege\nege J$ and $\vdash_i Q\to Q^{\Em}$. Then, as \\  $\vdash_i(J\to Q)\to (\nege\nege J\to \nege\nege Q)$ holds by Lemma \ref{I2}(2),
  $$(*)\;\;\;(J\to Q)\to (J^{\Em}\to \nege \nege Q).$$
  We also obtain $\vdash_i \nege \nege Q^{\Em}\to Q^{\Em}$ from  Corollary \ref{I5} (i). As $\vdash_i Q\to Q^{\Em}$ yields $\vdash_i \nege\nege Q\to \nege\nege Q^{\Em}$,
  we have $\nege\nege Q\to Q^{\Em}$, which yields \\  $\vdash_i(J\to Q) \to (J^{\Em}\to Q^{\Em})$ by $(*)$, thus $\vdash_i(J\to Q) \to (J\to Q)^{\Em}$. 
  \\[2ex]
  (ii): Since $\vdash_i\neg\perp$ we have $\vdash_i\nege\neg\perp\to\Em$.
  
  Now suppose that $\vdash_i\nege\neg R_j\to R^{\Em}_j$ holds for $j\in\{0,1\}$. According to Lemma \ref{I2}(3) we have 
  $$\vdash_i\nege \neg(R_0\wedge R_1)\to \nege \neg R_0\,\wedge\,\nege\neg R_1$$
  and thus $\vdash_i\nege\neg (R_0\wedge R_1)\to R_0^{\Em}\,\wedge\,R_1^{\Em}$, i.e., $\vdash_i\nege\neg (R_0\wedge R_1)\to (R_0\,\wedge\,R_1)^{\Em}$.
  
  We also have $\vdash_i\nege\neg (R_0\vee R_1)\to \nege\nege(\nege\neg R_0\,\vee\,\nege\neg R_1)$ by Lemma \ref{I2}(5), and hence 
  $$\vdash_i\nege\neg(R_0\vee R_1)\to \nege \nege (R_0^{\Em}\vee R_1^{\Em})$$
  i.e., $\vdash_i\nege\neg(R_0\vee R_1)\to  (R_0\vee R_1)^{\Em}$.
  
  Next assume that $\vdash_i\nege\neg R(a)\to R(a)^{\Em}$. By Lemma \ref{I2}(9) we have $\vdash_i\nege \neg \forall x R(x)\to \forall x \nege \neg R(x)$.
  Therefore, $\vdash_i \nege \neg \forall x R(x)\to \forall x R(x)^{\Em}$, i.e., $\vdash_i \nege \neg \forall x R(x)\to (\forall x R(x))^{\Em}$.
  
  Finally suppose that $\vdash_i J^{\Em}\to \nege\nege J$ and $\vdash_i\nege\neg R\to R^{\Em}$. Then, 
  $$\vdash_i (\nege\nege J\to \nege\neg R) \to(J^{\Em}\to R^{\Em})$$ and hence, by Lemma \ref{I2}(7),
   $\vdash_i \nege\neg( J\to R) \to(J^{\Em}\to R^{\Em})$, i.e., \\  $\vdash_i \nege\neg( J\to R) \to(J\to R)^{\Em}$.
   \\[2ex]
   (iii): We have $\vdash_i\perp^{\Em}\to \nege\nege \perp$ and $\vdash_i A^{\Em}\to \nege\nege A$ for atomic $A$ since $A^{\Em}\equiv \nege\nege A$. 
   
   Now assume that $\vdash_i J_i^{\Em}\to \nege\nege J_i$ for $i\in\{0,1\}$. Then, \\ $\vdash_i(J_0\wedge J_1)^{\Em}\to \nege\nege (J_0\wedge J_1)$, thus
   $\vdash_i(J_0\wedge J_1)^{\Em} \to \nege\nege (J_0\wedge J_1)$ follows by Lemma \ref{I2}(4). 
   
   We also have $\vdash_i J_0^{\Em}\,\vee\,J_1^{\Em}\to \nege\nege J_0\,\vee\,\nege\nege J_1$ and hence 
   $$\vdash_i \nege\nege(J_0^{\Em}\,\vee\,J_1^{\Em})\to \nege\nege(\nege\nege J_0\,\vee\,\nege\nege J_1),$$
  from which $\vdash_i(J_0\vee J_1)^{\Em}\to \nege\nege(J_0\vee J_1)$ follows by Lemma \ref{I2}(6). 
  
  Assuming $\vdash_i J(a)^{\Em}\to \nege\nege J(a)$, we have  $\vdash_i \exists x J(x)^{\Em}\to \exists x\,\nege\nege  J(x)$, and therefore 
  $\vdash_i \nege \nege \exists x J(x)^{\Em}\to \nege\nege \exists x  J(x)$ by Lemma \ref{I2}(10), i.e.,  \\ $\vdash_i (\exists x J(x))^{\Em}\to \nege\nege \exists x  J(x)$.
  
  Finally, assume $\vdash_i\nege \neg R\to R^{\Em}$ and $\vdash_i J^{\Em}\to \nege\nege J$. Then, $$\vdash_i(R^{\Em}\to J^{\Em}) \to \nege \neg R \to \nege\nege J,$$
  and hence, by Lemma \ref{I2}(8), $\vdash_i(R^{\Em}\to J^{\Em})\to \nege \nege (R\to J)$, i.e., \\ $\vdash_i(R\to J)^{\Em}\to \nege \nege (R\to J)$. 
  \qed

  We will make use of substitutions for variables and for the placeholder $\Em$.
  In the sequent calculi \`a la Gentzen and Takeuti \cite{takeuti87} and in the Sch\"utte calculi \cite{sch60,sch77} one distinguishes syntactically between free $a,b,c,\ldots$ and bound $x,y,z,\ldots$ variables.
  As terms can contain only free variables there will never be a problem of substitutability of terms for variables.
  
  We use $\varphi\subs at$ for the result of replacing every occurrence of the free variable $a$ in $\varphi$ by the term $t$. Similarly, for a sequent $\Gamma\Rightarrow \Delta$ and a derivation $\mathcal D$  we
  use $\Gamma\subs a t \Rightarrow \Delta \subs a t$ and $\mathcal{D}\subs a t$, respectively,  for the result of replacing every occurrence of the free variable $a$ by the term $t$.
  Note, however, that while $\varphi\subs t a$ will be a formula, too, $\mathcal{D}\subs a t$ may no longer be a derivation.
  
  In a similar vein, for a propositional constant  $\Em$ we convey the result of replacing each of its occurrences in a formula, sequent and derivation via  $\varphi\subs {\Em}{\psi}$,
   $\Gamma\subs {\Em}{\psi} \Rightarrow \Delta \subs{\Em}{\psi}$ and $\mathcal{D}\subs{\Em}{\psi}$, respectively. 
However, we have to add a caveat  here.   The formula formation rules in calculi that have different symbols for free and bound variables (such as \cite{takeuti87}) allow to go from a formula $\varphi(a)$ to $\forall x\varphi(x)$ only if $x$ does not already occur in $\varphi(a)$. Thus, if one substitutes a formula $\psi$ for $\Em$ in a formula the resulting syntactic object may no longer be a formula. As a result, we tacitly require that before substitutions are made, bound variables in $\psi$ have to replaced by ones that avoid this clash. 

Now, by obeying this additional requirement,   $\varphi\subs {\Em}{\psi}$ will be again a formula.  However, $\mathcal{D}\subs{\Em}{\psi}$ may no longer be a derivation 
  as some eigenvariable conditions may have become violated in the process. So we have to take care of that as well.

   Two substitutability results will be useful.
  
  \begin{lem}\label{subs1} Let $T$ be a theory. If $\mathcal{D}(a)$ is a $T$-derivation of $\Gamma\Rightarrow \Delta$ and $c$ is a variable that doesn't occur in $\mathcal{D}$ then
   $\mathcal{D}\subs a c$ is a
 $T$-derivation of $\Gamma\subs ac\Rightarrow \Delta\subs ac$.  \end{lem}
  \prf This is obvious as $c$ is a completely new variable as fas as $\mathcal D$ is concerned, so no eigenvarible conditions are affected by this substitution. Formally one proves this  by induction on the number of inferences of $\mathcal{D}$ (see \cite[CH. 1, Lemma 2.10]{takeuti87}\qed

   \begin{prop}\label{subs2} Let $T$ be a theory whose axioms do not contain $\Em$. If $\mathcal{D}$ is a $T$-derivation of $\Gamma\Rightarrow \Delta$  and $\psi$ is an arbitrary formula
   then there is a $T$-derivation  $\mathcal{D}'$  
   of $\Gamma\subs{\Em}{\psi}\Rightarrow \Delta\subs{\Em}{\psi}$.  \end{prop}
  \prf Use induction on the number of inferences of $\mathcal{D}$.
  The only kinds of inferences we need to look at are $\forall:\mbox{right}$ and $\exists:\mbox{left}$. So suppose the last inference of $\mathcal D$ was $\forall:\mbox{right}$ with premiss
   $\Gamma\Rightarrow \Delta', \varphi$ and conclusion $\Gamma\Rightarrow \Delta_0, \forall x\,\varphi\subs ax$, where $\Delta$ is $ \Delta_0, \forall x\,\varphi\subs ax$ and $a$ does not
  occur in $\Gamma\Rightarrow\Delta$.
  Let $\mathcal{D}_0$ be the immediate subderivation of $\mathcal D$ with end sequent $\Gamma\Rightarrow \Delta_0, \varphi$.
  Let $c$ be a free variable that neither occurs in $\mathcal D$ nor in $\psi$. By Lemma \ref{subs1},  $\mathcal{D}_1 :=\mathcal{D}_0\subs ac$ is derivation, too.
  Note that $\mathcal{D}_1$ is a derivation of $$\Gamma\Rightarrow \Delta_0,\varphi\subs ac$$ owing to the eigenvariable condition satisfied by $a$.

  Since $\mathcal{D}_1$ has fewer inferences than $\mathcal D$ we can apply the induction hypothesis to arrive at a derivation $\mathcal{D}_2$ of
  $$\Gamma\subs {\Em}{\psi}\Rightarrow \Delta_0\subs {\Em}{\psi}, (\varphi\subs ac)\subs{\Em}{\psi}.$$

  As $c$ does not occur in $\Gamma\subs {\Em}{\psi}\Rightarrow \Delta_0\subs {\Em}{\psi}$ and $\psi$, we can apply an inference $\forall:\mbox{right}$
  to obtain a derivation $\mathcal{D}'$ 
  of $$\Gamma\subs {\Em}{\psi}\Rightarrow \Delta_0\subs {\Em}{\psi},  \forall x\,(((\varphi\subs ac)\subs {\Em}{\psi})\subs cx)$$

  which is the same as $\Gamma\subs{\Em}{\psi}\Rightarrow \Delta\subs{\Em}{\psi}$ since 
 $\forall x\,(((\varphi\subs ac)\subs {\Em}{\psi})\subs cx)\equiv  ( \forall x\,\varphi\subs ax)\subs{\Em}{\psi}$. 
  
   $\exists:\mbox{left}$ is dealt with in a similar fashion.
   
   We still haven't strictly shown that the proof length increases at most polynomially. This will be addressed in the next section.
   \qed


  \begin{thm}\label{GeometricTheorem} If $T$ is a geometric theory, i.e. the axioms of $T$ are geometric implications, and $\fa$ is a geometric implication, then
  
  $$T\vdash_c \fa\;\mbox{ yields }\; T\vdash_i\fa.$$
  Moreover, if $\mathcal D$ is a classical deduction of $\fa$ in $T$, then the size of the intuitionistic  deduction of $\fa$ in $T$ increases at most polynomially in the size of 
  $\mathcal D$.
  \end{thm}
  \prf Suppose $T\vdash_c\fa$.   By Corollary \ref{I5}(ii) we conclude that \begin{eqnarray}\label{haupt1} &&T^\Em\vdash_i\fa^\Em.\end{eqnarray}
   As $T\subseteq \IQ$ it follows from (\ref{haupt1}) and Proposition \ref{I7}(i) that
    \begin{eqnarray}\label{haupt2} &&T\vdash_i\fa^\Em.\end{eqnarray}
    Now $\fa$ is of the form $\forall x_1\ldots\forall x_r(\fb\to\fc)$ with $\fb,\fc\in \IQ\cap \IJJ$. Thus, by  Proposition \ref{I7}(i), we can conclude from (\ref{haupt1}) that
    
   \begin{eqnarray}\label{haupt3} &&T\vdash_i\fb\to \fc^\Em.\end{eqnarray}
  (\ref{haupt3}) in conjunction with Proposition \ref{I7}(iii) yields
   \begin{eqnarray}\label{haupt4} &&T\vdash_i\fb\to \nege\nege \fc.\end{eqnarray}
   Now, since $\Em$ is just a placeholder we may substitute $\fc$ for $\Em$ everywhere in the derivation of (\ref{haupt4}) by Proposition \ref{subs2}, yielding 
   a derivation showing that
     \begin{eqnarray}\label{haupt5} &&T\vdash_i\fb\to ((\fc\to\fc)\to \fc).\end{eqnarray}
     As a result of (\ref{haupt5}) we have $T\vdash_i\fb\to \fc$ and hence $T\vdash_i\fa$, as desired.
     \qed
      From the proof of the foregoing theorem it's clear that conservativity obtains for a wider collection of theories than just geometric ones. 
   
   \begin{cor}  If $T$ is a theory whose axioms are in $\IQ$  and $\fa$ is a geometric implication, then
  
  $$T\vdash_c \fa\;\mbox{ yields }\; T\vdash_i\fa.$$
  Moreover, if $\mathcal D$ is a classical deduction of $\fa$ in $T$, then the size of the intuitionistic  deduction of $\fa$ in $T$ increases at most polynomially in the size of 
  $\mathcal D$.
  \end{cor}
  \prf This follows from the proof of Theorem \ref{GeometricTheorem}. \qed
     
     \section{Polynomial time bounds}
In view of the foregoing results, it might be rather obvious that the transformation      of a classical proof of a geometric implication in a geometric theory  into an intuitionistic proof can be carried out in polynomial time.  It might be in order, though, to be a bit more precise. The plan, however, is not to do this in detail but rather from a ``higher'' point of view.
It is a fact that the syntax of first-order logic can be recognized and manipulated by polynomial time algorithms in Buss'   theory $S^1_2$. One place where the arithmetization of metamathematics for the sequent calculus is carried out in detail is  \cite[Ch. 7]{buss86}. Among other things, we require functions for the arithmetization of substitutions in Lemma \ref{subs1} and Proposition \ref{subs2}. They give rise to $\Sigma^b_1$-defined functions of $S^1_2$ (see \cite[p. 130]{buss86}, where this is carried for substitution of a term into a formula). Moreover, all 
 $\Sigma^b_1$-definable functions of $S^1_2$ are polynomial time computable functions (see \cite[Corollary 8]{buss86}).
As all manipulation of proofs in this paper can be carried out by $\Sigma^b_1$-definable functions of $S^1_2$ we have achieved our goal. 


\section{$\infty$-geometric theories} 
\subsection{The infinite geometric case}\label{infgep}
Much more powerful  notions of geometricity are available in infinitary logics.  $\mathcal{L}_{\infty\omega}$-logic allows for the formation of  infinite disjunctions $\bigvee \Phi$ and conjunctions $\bigwedge\Phi$, where $\Phi$ is an arbitrary set of (infinitary) formulae. 
 In this richer syntax a formula is said to be an {\em $\infty$-positive} formula,  if it can be generated from atoms and $\perp$ via
 $\vee,\wedge,\exists$ and $\bigvee$. More precisely, the latter means that the 
 infinite disjunction $\bigvee \Phi$ is an $\infty$-positive formula whenever  $\Phi$ is a set of $\infty$-positive formula.
 
 The $\infty$-geometric implications are obtained in the same way from the $\infty$-positive formula as the geometric implications are obtained from the positive formulas.
An $\infty$-geometric theory is one whose axioms are  $\infty$-geometric implications.

Examples of such theories are the theories of {\em flat modules} over a ring (see \cite{wraith}), {\em torsion groups}, fields prime characteristic, archimedean ordered fields and connected graphs. 
Even Peano arithmetic has an $\infty$-geometric axiomatization (see \cite[2.4]{rathjen2016}).

$\infty$-geometric classical theories are also conservative over their intuitionistic version for $\infty$-geometric formulas. This can be proved in Constructive Zermelo-Fraenkel set theory
$\mathbf{CZF}$ (see \cite[Theorem 7.9]{rathjen2016}) via cut elimination for $\mathcal{L}_{\infty\omega}$. The techniques of this paper can also be extended to the $\mathcal{L}_{\infty\omega}$ context. As a result one can prove conservativity already in a much weaker fragment of $\mathbf{CZF}$, namely intuitionistic Kripke-Platek set theory with elementhood induction restricted to $\Sigma$-formulas, $\mathbf{IKP}^r$.

As the theories $\CZF$ and $\mathbf{IKP}^r$ allow for witness extraction from proofs of existential statements (see \cite[6.1]{rathjen2012}, \cite[2.35]{cook-rathjen2016}, \cite{rathjen2021}) this offers the exciting prospect  of  extracting  bounds from proofs in classical
$\infty$-geometric theories. However, note that for this the classical proof must exist as an object in the constructive background theory. So it's not enough to know of its existence 
by appealing to principles such as the axiom of choice or Zorn's lemma.

\paragraph{Acknowledgement}
The author acknowledges support via the John Templeton Foundation grant ``A new dawn of intuitionism''.


\end{document}

%% file: Barr_Remarks_macros.tex
\newcommand{\comment}[1]{}

\newtheorem{thm}{Theorem:}[section]

\newtheorem{cor}[thm]{Corollary:}
\newtheorem{prop}[thm]{Proposition:}

\newtheorem{lem}[thm]{Lemma:}
\newtheorem{deff}[thm]{Definition:}

\newtheorem{examples}[thm]{Examples:}

\newcommand{\CZF}{{\mathbf{CZF}}}

\newcommand{\prf}{{\bf Proof: }}

\newcommand{\beqs}{\begin{eqnarray*}}
\newcommand{\eeqs}{\end{eqnarray*}}

\newcommand{\beq}{\begin{eqnarray}}
\newcommand{\eeq}{\end{eqnarray}}

\newcommand{\suc}{{\mathrm{suc}}}

\newcommand{\PAs}[1]
{#1^+}

\newcommand{\bes}{\begin{eqnarray*}}
\newcommand{\ees}{\end{eqnarray*}}

\newcommand{\bgf}{\begin{frame}}
\newcommand{\efr}{\end{frame}}

\newcommand{\ee}{\end{itemize}\end{frame}}

\newcommand{\subs}[2]{({#1}/{#2})}